\documentclass[12pt]{article}
\usepackage{latexsym,amsmath}
\usepackage{times}
\usepackage{url}
\usepackage{amssymb}

\newtheorem{thm}{Theorem}[section]

\newtheorem{lemma}[thm]{Lemma}

\newtheorem{cor}[thm]{Corollary}

\newcommand{\beq}[1]{\begin{equation}\label{#1}}
\newcommand{\enq}[0]{\end{equation}}

\newcommand{\qed}[0]{{\hspace*{\fill}\mbox{$\Box$}}}

\newcommand{\ant}{{\operatorname{\rm int}\,}}

\newcommand{\ul}[0]{\underline}

\newcommand{\cA}[0]{{\cal A}}

\newcommand{\cE}[0]{{\cal E}}
\newcommand{\cI}[0]{{\cal I}}
\newcommand{\cJ}[0]{{\cal J}}
\newcommand{\cM}[0]{{\cal M}}
\newcommand{\cO}[0]{{\cal O}}

\newcommand{\cS}[0]{{\cal S}}
\newcommand{\cT}[0]{{\cal T}}
\newcommand{\cU}[0]{{\cal U}}
\newcommand{\cV}[0]{{\cal V}}
\newcommand{\cW}[0]{{\cal W}}

\newcommand{\ga}[0]{\alpha}
\newcommand{\gb}[0]{\beta}
\newcommand{\gre}[0]{\epsilon}

\newcommand{\gD}[0]{\Delta}
\newcommand{\grg}[0]{\gamma}

\newcommand{\gl}[0]{\lambda}
\newcommand{\go}[0]{\omega}
\newcommand{\gO}[0]{\Omega}
\newcommand{\gS}[0]{\Sigma}
\newcommand{\gs}[0]{\sigma}

\begin{document}

\renewcommand{\thefootnote}{\fnsymbol{footnote}}

\title{Sampling independent sets in the discrete torus}

\date{Department of Mathematics\\University of
Notre Dame\\South Bend, IN 46556\\\url{dgalvin1@nd.edu}}

\author{David Galvin\thanks{Work partially
supported by the National Science Foundation under agreement
DMS-0111298 while the author was a member of the Institute for
Advanced Study in Princeton, New Jersey.}}

\maketitle

\begin{abstract}
The even discrete torus is the graph $T_{L,d}$ on vertex set
$\{0,\ldots,L-1\}^d$ (with $L$ even) in which two vertices are
adjacent if they differ on exactly one coordinate and differ by $1
~(\mbox{mod } L)$ on that coordinate. The {\em hard-core measure
with activity $\gl$} on $T_{L,d}$ is the probability distribution
$\pi_\gl$ on the independent sets (sets of vertices spanning no
edges) of $T_{L,d}$ in which an independent set $I$ is chosen with
probability proportional to $\gl^{|I|}$. This distribution occurs
naturally in problems from statistical physics and the study of
communication networks.

We study Glauber dynamics, a single-site update Markov chain on the
set of independent sets of $T_{L,d}$ whose stationary distribution
is $\pi_\gl$. We show that for $\gl = \omega(d^{-1/4}\log^{3/4}d)$
and $d$ sufficiently large the convergence to stationarity is
(essentially) exponentially slow in $L^{d-1}$. This improves a
result of Borgs {\em et al.}, who had shown slow mixing of Glauber
dynamics for $\gl$ growing exponentially with $d$.

Our proof, which extends to $\rho$-local chains (chains which alter
the state of at most a proportion $\rho$ of the vertices in each
step) for suitable $\rho$, closely follows the conductance argument
of Borgs {\em et al.}, adding to it some combinatorial enumeration
methods that are modifications of those used by Galvin and Kahn to
show that the hard-core model with parameter $\gl$ on the integer
lattice ${\mathbb Z}^d$ exhibits phase coexistence for $\gl =
\omega(d^{-1/4}\log^{3/4}d)$.

The discrete even torus is a bipartite graph, with partition classes
$\cE$ (consisting of those vertices the sum of whose coordinates is
even) and $\cO$. Our result can be expressed combinatorially as the
statement that for each sufficiently large $\gl$, there is a
$\rho(\gl)>0$ such that if $I$ is an independent set chosen
according to $\pi_\gl$, then the probability that $||I \cap \cE|-|I
\cap \cO||$ is at most $\rho(\gl)L^d$ is exponentially small in
$L^{d-1}$. In particular, we obtain the combinatorial result that
for all $\varepsilon
>0$ the probability that a uniformly chosen independent set from
$T_{L,d}$ satisfies $||I \cap \cE|-|I \cap \cO|| \leq (.25 -
\varepsilon)L^d$ is exponentially small in $L^{d-1}$.

\end{abstract}

\section{Introduction and statement of the result}
\label{sec-intro}

Let $\gS=(V,E)$ be a simple, loopless, finite graph on vertex set
$V$ and edge set $E$. (For graph theory basics, see {\em e.g.}
\cite{Bollobas}, \cite{Diestel}.) Write $\cI(\gS)$ for the set of
independent sets (sets of vertices spanning no edges) in $V$. For
$\gl > 0$ we define the {\em hard-core probability measure with
activity $\gl$} on $\cI(\gS)$ by
$$
\pi_{\gl}(\{I\}) = \frac{\gl^{|I|}}{Z_{\gl}(\gS)} ~~~\mbox{for
$I\in \cI(\gS)$}
$$
where $Z_{\gl}(\gS) = \sum_{I\in \cI}\gl^{|I|}$ is the appropriate
normalizing constant or {\em partition function}. Note that $\pi_1$
is uniform measure on $\cI(\gS)$.

The hard-core measure originally arose in statistical physics (see
{\em e.g.} \cite{Dobrushin,BergSteif}) where it serves as a model of
a gas with particles of non-negligible size. The vertices of $\gS$
we think of as sites that may or may not be occupied by particles;
the rule of occupation is that adjacent sites may not be
simultaneously occupied. In this context the activity $\gl$ measures
the likelihood of a site being occupied.

The measure also has a natural interpretation in the context of
multicast communications networks (see {\em e.g.} \cite{Kelly}).
Here the vertices of $\gS$ are thought of as locations from which
calls can be made; when a call is made, the call location is
connected to all its neighbours, and throughout its duration, no
call may be placed from any of the neighbours. Thus at any given
time, the set of locations from which calls are being made is an
independent set in $\gS$. If calls are attempted independently at
each vertex as a Poisson process of rate $\lambda$ and have
independent exponential mean $1$ lengths, then the process has
stationary distribution $\pi_\gl$.

\medskip

Unless $L$ and $d$ are small, it is unfeasible to explicitly compute
the partition function $Z_\gl$ and the distribution $\pi_\gl$. It is
therefore of great interest to understand the effectiveness of
algorithms which approximate $Z_\gl$ and/or $\pi_\gl$. In this paper
we study {\em Glauber dynamics}, a Monte Carlo Markov chain (MCMC)
which simulates $\pi_\gl$. MCMC's occur frequently in computer
science in algorithms designed to sample from or estimate the size
of large combinatorially defined structures; they are also used in
statistical physics and the study of networks to help understand the
behavior of models of physical systems and networks in equilibrium.
Glauber dynamics is the single-site update Markov chain
$\cM_{\gl}=\cM_{\gl}(\gS)$ on state space $\cI(\gS)$ with transition
probabilities $P_{\gl}(I,J), ~I,J \in \cI(\gS),$ given by
$$
P_{\gl}(I,J) = \left\{
            \begin{array}{ll}
               0 & \mbox{ if $|I \bigtriangleup J| > 1$} \\
               \frac{1}{|V|}\frac{\gl}{1+\gl}& \mbox{ if $|I \bigtriangleup
J| = 1, ~I \subseteq
               J$} \\
               \frac{1}{|V|}\frac{1}{1+\gl}& \mbox{ if $|I \bigtriangleup J|
= 1, ~J \subseteq
               I$} \\
               1 - \sum_{I \neq J' \in \cI(\gS)} P_{\gl}(I,J') & \mbox{ if
$I=J$.}
            \end{array}
         \right.
$$
We may think of $\cM_{\gl}$ dynamically as follows. From an
independent set $I$, choose a vertex $v$ uniformly from $V$. Then
add $v$ to $I$ with probability proportional to $\gl$, and remove it
with probability proportional to $1$; that is, set
$$
I' = \left\{
         \begin{array}{ll}
             I \cup \{v\} & \mbox{ with probability $\frac{\gl}{1+\gl}$} \\
             I \setminus \{v\} & \mbox{ with probability
             $\frac{1}{1+\gl}$}.
         \end{array}
     \right.
$$
Finally, move to $I'$ if $I'$ is an independent set, and stay at
$I$ otherwise.

It is readily checked that $\cM_{\gl}$ is an ergodic Markov chain
with (unique) stationary distribution $\pi_{\gl}$. A natural and
important question to ask about $\cM_{\gl}$ is how quickly it
converges to its stationary distribution. It is traditional to
define the {\em mixing time} $\tau_{{\cal M}_{\gl}(\gS)}$ of ${\cal
M}_{\gl}(\gS)$ to be
$$
\tau_{\cM_\gl(\gS)}=\max_{I \in \cI(\gS)} \min
\left\{t_0:\frac{1}{2}\sum_{J \in \cI(\gS)}|P^t(I,J)-\pi_\gl(J)|
\leq \frac{1}{e} ~~~ \forall t>t_0\right\},
$$
where $P^{t}(I, \cdot)$ is the distribution of the chain at time
$t$, given that it started in state $I$. The mixing time of
$\cM_\gl$ captures the speed at which the chain converges to its
stationary distribution: for every $\gre
>0$, in order to get a sample from $\cI(\gS)$ which is within $\gre$
of $\pi_\gl$ (in variation distance), it is necessary and sufficient
to run the chain from some arbitrarily chosen distribution for some
multiple (depending on $\gre$) of the mixing time. For surveys of
issues related to the mixing time of a Markov chain, see {\em e.g.}
\cite{MontenegroTetali,RandallMixing}.

\medskip

Here we study $\tau_{\cM_\gl(T_{L,d})}$, where $T_{L,d}$ is the even
discrete torus. This is the graph on vertex set $\{0,\ldots,L-1\}^d$
(with $L$ even) in which two strings are adjacent if they differ on
only one coordinate, and differ by $1 ~(\mbox{mod } L)$ on that
coordinate. For $L\geq 4$ this is a $2d$-regular bipartite graph
with unique bipartition $\cE \cup \cO$ where ${\cal E}$ is the set
of even vertices of $T_{L,d}$ (those strings the sum of whose
coordinates is even) and ${\cal O}$ is the set of odd vertices.

Much work has been done on the question of bounding $\tau_{\cM_\gl}$
above for various classes of graphs. The most general results
available to date are due to Luby and Vigoda \cite{LubyVigoda} and
Dyer and Greenhill \cite{DyerGreenhill}, who have shown that for any
graph $\Sigma$ with maximum degree $\Delta$,
$\tau_{\cM_\gl(\Sigma)}$ is a polynomial in $|V(\Sigma)|$ whenever
$\gl < 2/(\Delta-2)$, which implies that $\tau_{\cM_\gl(T_{L,d})}$
is a polynomial in $L^d$ whenever $\gl < 1/(d-1)$. More recently,
Weitz \cite{Weitz} has improved this general bound in the case of
graphs with sub-exponential growth, and in particular has shown that
$\tau_{\cM_\gl(T_{L,d})}$ is a polynomial in $L^d$ whenever $\gl
\leq (2d-1)^{2d-1}/(2d-2)^{2d} \approx e/2d$.

Recently, attention has been given to the question of regimes of
{\em inefficiency} of Glauber and other dynamics. Dyer, Frieze and
Jerrum \cite{DyerFriezeJerrum} considered the case $\gl=1$ and
showed that for each $\gD\geq 6$ a random (uniform) $\gD$-regular,
$n$-vertex bipartite $\gS$ almost surely (with probability tending
to $1$ as $n$ tends to infinity) satisfies $\tau_{\cM_1}(\gS) \geq
2^{\gamma n}$ for some absolute constant $\gamma >0$. The first
result in this vein that applied specifically to $T_{L,d}$ was due
to Borgs {\em et al.} \cite{BorgsChayesFriezeKimTetaliVigodaVu}, who
used a conductance argument to obtain the following.
\begin{thm} \label{thm-from.million}
There is $c(d)>0$ (independent of $L$) such that for $\gl$
sufficiently large and all even $L \geq 4$,
$$
\tau_{\cM_\gl(T_{L,d})} > \exp\left\{\frac{c(d)L^{d-1}}{\log ^2
L}\right\}.
$$
\end{thm}

An examination of \cite{BorgsChayesFriezeKimTetaliVigodaVu} reveals
that ``sufficiently large'' may be quantified as $\gl > c^d$ for a
suitable constant $c>1$. One motivation for
\cite{BorgsChayesFriezeKimTetaliVigodaVu} was to show that for
values of $\gl$ for which the hard-core model on the integer lattice
${\mathbb Z}^d$ exhibits multiple Gibbs phases (to be explained
below), the mixing of the Glauber dynamics on $T_{L,d}$ should be
slow. Dobrushin \cite{Dobrushin} showed that as long as $\gl$ is
sufficiently large, there are indeed multiple Gibbs phases in the
hard-core model. Specifically, write $\cE$ and $\cO$ for the sets of
even and odd vertices of ${\mathbb Z}^d$ (defined in the obvious
way). Equip ${\mathbb Z}^d$ with the usual nearest neighbour
adjacency and set
$$
\Lambda_L =[-L,L]^d ~~~\mbox{and}~~~ \partial \Lambda_L =[-L,L]^d
\setminus [-(L-1),L-1]^d.
$$
For $\lambda>0$, choose ${\mathbb I}$ from ${\cal I}(\Lambda_L)$
with $\Pr({\mathbb I}=I) \propto \lambda^{|I|}$. Dobrushin showed
that for $\lambda$ large
\begin{equation} \label{thm-Dob}
\lim_{L\rightarrow\infty}{\mathbb P}\left(\vec{0}\in{\mathbb
I}~|~{\mathbb I}\supseteq
\partial \Lambda_L\cap {\cal E}\right)~>
\lim_{L\rightarrow\infty}{\mathbb P}\left(\vec{0}\in{\mathbb
I}~|~{\mathbb I}\supseteq
\partial \Lambda_L\cap {\cal O}\right)
\end{equation}
where $\vec{0}=(0,\ldots, 0)$. Thus, roughly speaking, the influence
of the boundary on behavior at the origin persists as the boundary
recedes. Informally, this suggests that for $\gl$ large, the typical
independent set chosen from $T_{L,d}$ according to the hard-core
measure is either predominantly odd or predominantly even, and so
there is a highly unlikely bottleneck set of balanced independent
sets separating the predominantly odd sets from the predominantly
even ones. It is the existence of this bottleneck that should cause
the mixing of the Glauber dynamics chain to be slow. No explicit
bound is given in \cite{Dobrushin}, but several researchers report
that Dobrushin's argument works for $\gl>c^d$ for a suitable
constant $c>1$. A key tool in the proof of Theorem
\ref{thm-from.million} is an appeal to a (suitable generalization)
of a lemma of Dobrushin from \cite{Dobrushin2}, and our main lemma,
Lemma \ref{lem-volume.bounds.from.gk}, is of a similar flavour.

In light of a recent result of Galvin and Kahn \cite{GalvinKahn}, it
is tempting to believe that slow mixing on $T_{L,d}$ should hold for
smaller values of $\gl$; even for values of $\gl$ tending to $0$ as
$d$ grows. The main result of \cite{GalvinKahn} is that the
hard-core model on ${\mathbb Z}^d$ exhibits multiple Gibbs phases
for $\gl=\omega(d^{-1/4}\log^{3/4}d)$. Specifically, Galvin and Kahn
show that for $\gl \geq cd^{-1/4}\log^{3/4}d$ for sufficiently large
$c$, (\ref{thm-Dob}) holds.

In \cite{GalvinTetali}, some progress was made towards establishing
slow mixing on $T_{L,d}$ for small $\gl$. Let $Q_d$ be the usual
discrete hypercube (the graph on $\{0,1\}^d$ in which two strings
are adjacent if they differ on exactly one coordinate). Note that
$T_{2,d}$ is isomorphic to $Q_d$. A corollary of the main result of
\cite{GalvinTetali} is that for $\gl = \omega(d^{-1/4}\log^{3/2}d)$,
$$
\tau_{\cM_\gl(Q_d)} \geq \exp
\left\{\Omega\left(\frac{2^d}{d^2}\right)\right\}.
$$

In the present paper, using different methods, we show that for
$d$ sufficiently large Glauber dynamics does indeed mix slowly on
$T_{L,d}$ for all even $L \geq 4$ for some small values of
$\lambda$.
\begin{thm} \label{thm-main}
There are constants $c, d_0>0$ for which the following holds. For
\begin{equation} \label{inq-bound.on.gl}
\gl \geq cd^{-1/4}\log^{3/4}d,
\end{equation}
$d\geq d_0$ and $L\geq 4$ even, the Glauber dynamics chain $\cM_\gl$
on $\cI(T_{L,d})$ satisfies
$$
\tau_{\cM_\gl(T_{L,d})} \geq \exp \left\{\frac{L^{d-1}}{d^4\log^2
L}\right\}.
$$
\end{thm}

\medskip

Our techniques actually apply to the class of $\rho$-local chains
(considered in \cite{BorgsChayesFriezeKimTetaliVigodaVu} and also in
\cite{DyerFriezeJerrum}, where the terminology {\em
$\rho|V|$-cautious} is employed) for suitable $\rho$. A Markov chain
$\cM$ on state space $\cI$ is {\em $\rho$-local} if in each step of
the chain the states of at most $\rho|V|$ vertices are changed; that
is, if
$$
P_\cM(I_1,I_2) \neq 0 \Rightarrow |I_1 \triangle I_2| \leq \rho|V|.
$$
Our main theorem is the following.

\begin{thm} \label{thm-main1.5}
There are constants $c, d_0>0$ for which the following holds. For
$\gl$ satisfying (\ref{inq-bound.on.gl}), $d\geq d_0$, $L\geq 4$
even and $\rho$ satisfying
\begin{equation} \label{inq-condition.on.alpha.rho.1}
\rho + \frac{1}{2d^{1/2}} \leq \frac{\gl}{1+\gl}
\end{equation}
and
\begin{equation} \label{inq-condition.on.alpha.rho}
H\left(\frac{1}{2d^{1/2}}\right)+H\left(\rho +
\frac{1}{2d^{1/2}}\right)+\left(\frac{1}{d^{1/2}}+\rho\right)\log_2
\gl + \frac{10}{d^4L \log^2 L} \leq \log_2(1+\gl)
\end{equation}
(where $H(\ga)=-\ga\log_2\ga-(1-\ga)\log_2(1-\ga)$ is the usual
binary entropy function), if $\cM$ is an ergodic $\rho$-local Markov
chain on state space $\cI(T_{L,d})$ with stationary distribution
$\pi_\gl$ then
$$
\tau_{\cM(T_{L,d})} \geq \exp\left\{\frac{L^{d-1}}{d^4\log ^2
L}\right\}.
$$
\end{thm}

With $\rho=L^{-d}$, (\ref{inq-condition.on.alpha.rho}) is satisfied
for all $\gl$ satisfying (\ref{inq-bound.on.gl}) (for sufficiently
large $d$). An $L^{-d}$-local chain is a single-site update chain
and so Theorem \ref{thm-main} is a corollary of Theorem
\ref{thm-main1.5}. Taking $\gl=1$ we may satisfy
(\ref{inq-condition.on.alpha.rho}) with $\rho$ any constant less
than $1/2$ by taking $d$ large enough (as a function of $\rho$). We
therefore obtain a further corollary of Theorem \ref{thm-main1.5}.

\begin{cor} \label{cor-gl.equals.1}
Fix $\rho < 1/2$. There is a constant $d_0=d_0(\rho)>0$ for which
the following holds. For $L\geq 4$ even and $d \geq d_0$, if $\cM$
is an ergodic $\rho$-local Markov chain on state space
$\cI(T_{L,d})$ with uniform stationary distribution then
$$
\tau_{\cM(T_{L,d})} \geq \exp\left\{\frac{L^{d-1}}{d^4\log ^2
L}\right\}.
$$
\end{cor}

\medskip

We prove Theorem \ref{thm-main1.5} via a well-known conductance
argument (introduced in \cite{JerrumSinclair}). A particularly
useful form of the argument was given by Dyer, Frieze and Jerrum
\cite{DyerFriezeJerrum}. Let $\cM$ be an ergodic Markov chain on
state space $\gO$ with transition probabilities $P$ and stationary
distribution $\pi$. Let $A \subseteq \gO$ and $M \subseteq \gO
\setminus A$ satisfy $\pi(A) \leq 1/2$ and
$$
\go_1 \in A, \go_2 \in \gO \setminus (A \cup M) \Rightarrow P(\go_1,
\go_2) =0.
$$
Then from \cite{DyerFriezeJerrum} we have
\begin{equation} \label{conductance}
\tau_\cM \geq \frac{\pi(A)}{8\pi(M)}.
\end{equation}
The intuition behind (\ref{conductance}) is that if we start the
chain at some state in $A$, then in order to mix, it must at some
point leave $A$ and so pass through $M$. The ratio of $\pi(A)$ to
$\pi(M)$ is a measure of how long the chain must run before it
transitions from $A$ to $M$. So we may think of $M$ as a bottleneck
set through which any run of the chain must pass in order to mix; if
the bottleneck has small measure, then the mixing time is high.

\medskip

Now let us return to the setup of Theorem \ref{thm-main1.5}. Set
$$
\cI_{b,\rho} = \cI_{b,\rho}(T_{L,d}) = \{I \in
\cI(T_{L,d}):\left||I \cap \cE|-|I \cap \cO|\right|\leq \rho
L^d/2\}
$$ ($\cI_{b,\rho}$ is the set
of {\em balanced} independent sets) and
$$
\cI_{\cE,\rho} = \cI_{\cE,\rho}(T_{L,d}) = \{I \in \cI(T_{L,d}):|I
\cap \cE|>|I \cap \cO|+\rho L^d/2\}.
$$
By symmetry, $\pi_\gl(\cI_{\cE,\rho})<1/2$. Notice that since $\cM$
changes the state of at most $\rho L^d$ vertices in each step, we
have that if $I_1 \in \cI_{\cE,\rho}$ and $I_2 \in \cI(T_{L,d})
\setminus (\cI_{\cE,\rho} \cup \cI_{b,\rho})$ then $P_\cM(I_1,I_2) =
0$. From (\ref{conductance}) we obtain
$$
\tau_\cM \geq \frac{\pi_\gl(\cI_{\cE,\rho})}{8\pi_\gl(\cI_{b,\rho})}
= \frac{1-\pi_\gl(\cI_{b,\rho})}{16\pi_\gl(\cI_{b,\rho})}.
$$
Theorem \ref{thm-main1.5} thus follows from the following theorem,
whose proof will be the main business of this paper.
\begin{thm} \label{thm-main2}
There are constants $c, d_0 > 0$ for which the following holds. For
$\gl$ satisfying (\ref{inq-bound.on.gl}), $d\geq d_0$, $L\geq 4$
even and $\rho$ satisfying (\ref{inq-condition.on.alpha.rho.1}) and
(\ref{inq-condition.on.alpha.rho}),
$$
\pi_\gl(\cI_{b,\rho}) \leq \exp\left\{-\frac{2L^{d-1}}{d^4\log ^2
L}\right\}.
$$
\end{thm}

Theorem \ref{thm-main2} is the statement that if an independent set
$I$ is chosen from $\cI(T_{L,d})$ according to the hard-core
distribution $\pi_\gl$, then, as long as $\gl$ is sufficiently
large, it is extremely unlikely that $I$ is balanced. In particular,
if we take $\gl=1$ we obtain the following appealing combinatorial
corollary.
\begin{cor} \label{cor-balanced}
Fix $\varepsilon > 0$. There is a constant $d_0=d_0(\varepsilon)>0$
for which the following holds. For $L\geq 4$ even and $d \geq d_0$,
if ${\mathbb I}$ is a uniformly chosen independent set from
$T_{L,d}$ then
$$
{\mathbb P}\left(\left||{\mathbb I} \cap \cE|-|{\mathbb I} \cap
\cO|\right| \leq (.25 - \varepsilon)L^d\right) \leq
\exp\left\{-\frac{2L^{d-1}}{d^4\log ^2 L}\right\}.
$$
\end{cor}

\section{Overview of the proof of Theorem \ref{thm-main2}}
\label{sec-overview}

Consider an independent set $I\in \cI(T_{L,d})$. Some regions of
$T_{L,d}$ consist predominantly of even vertices from $I$ together
with their neighbours (the {\em even-occupied} regions) and some
regions consist predominantly of odd vertices from $I$ with their
neighbours. These regions are separated by a collection of connected
unoccupied two-layer moats or {\em cutsets} $\grg$. In Section
\ref{subsec-cutsets} we follow
\cite{BorgsChayesFriezeKimTetaliVigodaVu} and describe a procedure
which selects a collection $\Gamma(I)$ of these $\grg$'s with the
properties that $i)$ the interiors of those $\grg \in \Gamma(I)$ are
mutually disjoint (the {\em interior} of $\grg$ is the smaller of
the two parts into which its deletion breaks a graph) and $ii)$
either the interiors of all $\grg \in \Gamma(I)$ are predominantly
even-occupied or they are all predominantly odd-occupied. We do this
in the setting of an arbitrary bipartite graph. We also point out
some properties of $\grg$ that are specific to the torus, including
an isoperimetric inequality that gives a lower bound on $|\grg|$
(the number of edges in $\grg$) in terms of the number of vertices
it encloses.

Our main technical result, Lemma \ref{lem-volume.bounds.from.gk}, is
the assertion that for each specification of cutset sizes $c_1,
\ldots, c_\ell$ and vertices $v_1, \ldots, v_\ell$, the probability
that an independent set $I$ has among its associated cutsets
$\Gamma(I)$ a collection $\grg_1, \ldots, \grg_\ell$ with
$|\grg_i|=c_i$ and with $v_i$ in the interior of $\grg_i$ is
exponentially small in the sum of the $c_i$'s. The case $\ell=1$ is
essentially contained in \cite{GalvinKahn}, and our generalization
draws heavily on that paper. It may be worthwhile to compare our
Lemma \ref{lem-volume.bounds.from.gk} with \cite[Lemma
6]{BorgsChayesFriezeKimTetaliVigodaVu} in which is obtained an
exponential bound on the probability of $I$ having a {\em
particular} collection of cutsets.

We use a Peierl's argument (see {\em e.g.} \cite{Grimmett}) to prove
Lemma \ref{lem-volume.bounds.from.gk}. For simplicity, we describe
the argument here for $\gl=1$. For fixed $c_1, \ldots, c_\ell$,
$v_1, \ldots, v_\ell$, let $\cI_{spec}$ be the collection of $I \in
\cI(T_{L,d})$ which have a collection of associated cutsets $\grg_1,
\ldots, \grg_\ell$ with $|\grg_i|=c_\ell$ and with $v_i$ in the
interior of $\grg_i$. For an $I \in \cI_{spec}$, fix one such
collection $\grg_1, \ldots, \grg_\ell$. By modifying $I$ carefully
in the interior of each $\gamma_i$ (specifically, by shifting $I$
one unit in a carefully chosen direction) we can identify a
collection of subsets $S_i$ of the vertices of $\grg_i$ with $|S_i|
= c_i/2d$ which can be added to the modified $I$, the resulting set
still being independent. (Here we exploit the fact that the cutset
can be thought of as two unoccupied layers separating the interior
from the exterior). By adding arbitrary subsets of each $S_i$ to the
modified $I$, we get a one-to-many map $\varphi$ from $\cI_{spec}$
to $\cI(T_{L,d})$ with $|\varphi(I)|$ exponential in the sum of the
$c_i$'s.

If the $\varphi(I)$'s would be disjoint for distinct $I$'s, we would
essentially be done, having shown that there are exponentially more
(in the sum of the $c_i$'s) independent sets than sets in
$\cI_{spec}$. To deal with the issue of overlaps between the
$\varphi(I)$'s, we define a flow $\nu:\cI_{spec}
\times\cI(T_{L,d})\rightarrow [0,\infty)$ supported on pairs $(I,J)$
with $J \in \varphi(I)$ in such a way that the flow out of every $I
\in \cI_{spec}$ is $1$. Any uniform bound we can obtain on the flow
into vertices of $\cI(T_{L,d})$ is then easily seen to be a bound on
$\pi_1(\cI_{spec})$.

We define the flow via a notion of approximation modified from
\cite{GalvinKahn}. To each cutset $\grg$ we associate a set
$A(\grg)$ which approximates the interior of $\grg$ in a precise
sense, in such a way that as we run over all possible $\grg$, the
total number of approximate sets used is small (and in particular,
much smaller than the total number of cutsets). There is a clear
trade-off here: the more precise the notion of approximation used,
the greater the number of approximate sets needed. Then for each $J
\in \cI(T_{L,d})$ and each collection of approximations $A_1,
\ldots, A_\ell$ we consider the set of those $I \in \cI_{spec}$ with
$J \in \varphi(I)$ and with $A_i$ the approximation to $\grg_i$. We
define the flow in such a way that if this set is large, then
$\nu(I,J)$ is small for each $I$ in the set. In this way we control
the flow into $J$ corresponding to each collection of approximations
$A_1, \ldots, A_\ell$; and since the total number of approximations
is small, we control the total flow into $J$.

In the language of statistical physics, there is a tradeoff between
{\em entropy} and {\em energy} that we need to control. Each $I \in
\cI_{spec}$ has high energy --- by the shift operation described
above, we can perturb it only slightly and map it to an
exponentially large collection of independent sets. But before
exploiting this fact to show that $\pi_1(\cI_{spec})$ is small, we
have to account for a high entropy term --- there are exponentially
many possible cutsets of size $c_i$ that could be associated with an
$I \in \cI_{spec}$. There are about $\exp\{\Omega(c_i \log d /d)\}$
cutsets of size $c_i$ (this count comes from \cite{LebowitzMazel}),
each one giving rise to about $\exp\{\Theta(c_i/d)\}$ independent
sets, so the entropy term exceeds the energy term and the Peierl's
argument cannot succeed. One way to overcome this problem is to
allow $\gl$ to grow exponentially with $d$, increasing the energy
term (the independent sets obtained from the shift are larger than
the pre-shifted sets, and so have greater weight) while not changing
the entropy term. This is the approach taken in
\cite{BorgsChayesFriezeKimTetaliVigodaVu}. Alternatively we could
try to salvage the argument for $\gl=1$ by somehow decreasing the
entropy term. This is where the idea of approximate cutsets comes
in. Instead of specifying a cutset $\grg_i$ by its $c_i$ edges, we
specify a connected collection of roughly $c_i/d^{3/2}$ vertices
nearby (in a sense to be made precise) to the cutset, from which a
good approximation to the cutset can be constructed in a specified
(algorithmic) way. Our entropy term drops to roughly
$\exp\{O(c_i\log d/d^{3/2})\}$, much lower than the energy term; so
much lower, in fact, that we can rescue the Peierl's argument for
values of $\gl$ tending to $0$ as $d$ grows. The bound
$\exp\{O(c_i\log d/d^{3/2})\}$ on the number of connected subsets of
$T_{L,d}$ of size $O(c_i/d^{3/2})$ is based on the fact that a
$\Delta$-regular graph has at most $2^{O(n\log \Delta)}$ connected
induced subgraphs of size $n$ passing through a fixed vertex.

The precise statement of Lemma \ref{lem-volume.bounds.from.gk}
appears in Section \ref{subsec-results.from.GalvinKahn} and the
proof appears in Section \ref{subsec-proofs}. It is here that the
precise notion of approximation used is given, together with the
verification that there is a $\nu$ that satisfies our diverse
requirements. We defer a more detailed discussion of the proof to
that section.

Given Lemma \ref{lem-volume.bounds.from.gk}, the proof of Theorem
\ref{thm-main2} is relatively straightforward. We begin by using a
naive count to observe that the total measure of those $I \in
\cI_{b,\rho}$ with $\min\{|I\cap\cE|,|I\cap\cO|\}\leq L^d/4d^{1/2}$
is exponentially small in $L^d$. This drives our specification of
$\rho$, which is chosen as large as possible so that the naive count
gives an exponentially small bound. This allows us in the sequel to
consider only those $I \in \cI(T_{L,d})$ with
$\min\{|I\cap\cE|,|I\cap\cO|\}> L^d/4d^{1/2}$. The naive count
consists of considering those subsets $X$ of $T_{L,d}$ with
$\min\{|X \cap \cE|,|X \cap \cO|\} \leq L^d/4d^{1/2}$ and $\max\{|X
\cap \cE|,|X \cap \cO|\} \leq L^d/4d^{1/2} + \rho L^d/2$, without
regard for whether $X \in \cI(T_{L,d})$.

It remains to consider the case where balanced $I$ satisfies
$\min\{|I\cap\cE|,|I\cap\cO|\}> L^d/4d^{1/2}$. In this case the
isoperimetric inequality in the torus allows us to conclude that
$\Gamma(I)$ contains a small subset of cutsets, all with similar
lengths, the sum of whose lengths is essentially $L^{d-1}$. We then
use Lemma \ref{lem-volume.bounds.from.gk} and a union bound to say
that the measure of the large balanced independent sets is at most
the product of a term that is exponentially small in $L^{d-1}$ (from
Lemma \ref{lem-volume.bounds.from.gk}), a term corresponding to the
choice of a fixed vertex in each of the interiors, and a term
corresponding to the choice of the collection of lengths. The second
term will be negligible because our special collection of contours
is small and the third will be negligible because the contours all
have similar lengths. The detailed proof appears in Section
\ref{subsec-proof_of_thm}.

\section{Proof of Theorem \ref{thm-main2}} \label{sec-proof.of.mn.thm}

\subsection{Cutsets} \label{subsec-cutsets}

We describe a way of associating with each $I \in \cI(T_{L,d})$ a
collection of minimal edge cutsets, following the approach of
\cite{BorgsChayesFriezeKimTetaliVigodaVu}. Much of the discussion is
valid for any bipartite graph, so we present it in that generality.

\medskip

Let $\gS=(V,E)$ be a connected bipartite graph on at least $3$
vertices with partition classes $\cE$ and $\cO$. For $X \subseteq
V$, write $\nabla(X)$ for the set of edges in $E$ which have one end
in $X$ and one end outside $X$; $\overline{X}$ for $V \setminus X$;
$\partial_{int}X$ for the set of vertices in $X$ which are adjacent
to something outside $X$; $\partial_{ext}X$ for the set of vertices
outside $X$ which are adjacent to something in $X$; $X^+$ for $X
\cup
\partial_{ext} X$; $X^\cE$ for $X \cap \cE$ and $X^\cO$ for $X \cap
\cO$. Further, for $x \in V$ set $\partial x=\partial_{ext}\{x\}$.
In what follows we abuse notation slightly, identifying sets of
vertices of $V$ and the subgraphs they induce.

For each $I \in \cI(\gS)$, each component $R$ of $(I^\cE)^+$ or
$(I^\cO)^+$ and each component $C$ of $\overline{R}$, set $\gamma =
\gamma_{RC}(I)=\nabla(C)$ and $W=W_{RC}(I)=\overline{C}$. Evidently
$C$ is connected, and $W$ consists of $R$, which is connected,
together with a number of other components of $\overline{R}$, each
of which is connected and joined to $R$, so $W$ is connected also.
It follows that $\grg$ is a minimal edge-cutset in $\gS$. Define the
{\em size} of $\grg$ to be $|\grg|=|\nabla C|~(=|\nabla(W)|)$.
Define $\ant \gamma$, the {\em interior} of $\gamma$, to be the
smaller of $C,W$ (if $|W|=|C|$, take $\ant \gamma = W$) and say that
$\grg$ is {\em enveloping} if $\ant \gamma = W$ (so that $R$, the
component that gives rise to $\grg$, is contained in the interior of
$\grg$). Say that $I$ is {\em even} (respectively, {\em odd}) if it
satisfies the following condition: for every component $R$ of
$(I^\cE)^+$ (respectively, $(I^\cO)^+$) there exists a component $C$
of $\overline{R}$ such that $\gamma_{RC}(I)$ is enveloping. Note
that there must be an unique such $C$ for each $R$ since the
components of $\overline{R}$ are disjoint and each one that gives
rise to an enveloping cutset must have more than $|V|/2$ vertices.

\begin{lemma} \label{lem-even_or_odd}
Each $I \in \cI(\gS)$ is either odd or even.
\end{lemma}

\noindent {\em Proof: }Suppose that $I$ is not even. Then there is a
component $R$ of $(\cI^\cE)^+$ such that for all components $C$ of
$\overline{R}$, $|C|<|V|/2$. Consider a component $R'$ of
$(\cI^\cO)^+$. It lies inside some component $C$ of $\overline{R}$,
so one of the components of $\overline{R'}$, say $C'$, contains
$\overline{C}$. Since $|\overline{C}|\geq |V|/2$ the cutset
$\gamma_{R'C'}(I)$ is enveloping. It follows that $I$ is odd. \qed

\begin{lemma} \label{lem-substantial}
For each even $I \in \cI(\gS)$ there is an associated collection
$\Gamma(I)$ of enveloping cutsets with mutually disjoint interiors
such that $I^\cE \subseteq \cup_{\grg \in \Gamma(I)} \ant \grg$.
\end{lemma}

\noindent {\em Proof: }Let $R_1, \ldots, R_m$ be the components of
$(I^\cE)^+$. For each $i$ there is one component, $C_i$ say, of
$\overline{R_i}$ such that $\gamma_i = \gamma_{R_iC_i}$ is
enveloping. We have $I^\cE \subseteq \cup_{i=1}^m \ant \grg_i$.

We claim that for each $i \neq j$ one of $\ant \grg_i \subseteq \ant
\grg_j$, $\ant \grg_i \supseteq \ant \grg_j$, $\ant \grg_i \cap \ant
\grg_j = \emptyset$ holds. To see this, we consider cases. If $R_j
\subseteq C'$ for some component $C'\neq C_i$ of $\overline{R_i}$
then $\ant \grg_j \subseteq C' \subseteq \ant
\grg_i~(=\overline{C_i})$. Otherwise, $R_j \subseteq C_i$. In this
case, either $C_j \subseteq C_i$ (so $\ant \grg_j \supseteq \ant
\grg_i$) or $C_j \supseteq \overline{C_i}$ (so $\ant \grg_j \cap
\ant \grg_i = \emptyset$). We may take
$$
\Gamma(I)=\{\grg_i~:~\mbox{for all $j \neq i$ either $\ant \grg_j
\subseteq \ant \grg_i$ or $\ant \grg_i \cap \ant \grg_j =
\emptyset$}\}.
$$
\qed

The following lemma identifies some key properties of $\grg \in
\Gamma(I)$ for even $I$. In the proof of Theorem \ref{thm-main2}
these properties only come into play through Lemma
\ref{lem-volume.bounds.from.gk}.

\begin{lemma} \label{lem-properties.of.contours_general}
For each even $I$ and $\grg \in \Gamma(I)$, we have the following.
\begin{equation} \label{contour.prop.2}
\partial_{int}W \subseteq \cO~~\mbox{and}~~\partial_{ext}W \subseteq
\cE;
\end{equation}
\begin{equation} \label{contour.prop.3}
\partial_{int}W \cap I =
\emptyset~~\mbox{and}~~\partial_{ext}W \cap I = \emptyset;
\end{equation}
\begin{equation} \label{contour.prop.4}
\forall x \in \partial_{int}W, ~\partial  x \cap W \cap I \neq
\emptyset
\end{equation}
and
\begin{equation} \label{contour.prop.5}
W^\cO = \partial_{ext} W^\cE~~\mbox{and}~~W^\cE=\left\{y \in
\cE:\partial y \subseteq W^\cO\right\}.
\end{equation}
\end{lemma}

\noindent {\em Proof: }We begin by noting that $\partial_{int} W
\subseteq
\partial_{int} R$ (specifically, $\partial_{int} W =
\partial_{int} R \cap
\partial_{ext}C = \partial_{ext}C$) and $\partial_{ext} W =
\partial_{int}C$. Since $\partial_{int} R \subseteq \cO$ and
$\partial_{int}C \subseteq \cE$, (\ref{contour.prop.2}) follows
immediately from these observations.

By construction, $R \cap \cO \cap I = \emptyset$, so $\partial_{int}
W \cap I = \emptyset$. If there is $x \in
\partial_{int}C \cap I$ then, since $x \in \cE$ and there is $y
\in R$ adjacent to $x$, we would have $x \in R$, a contradiction; so
$\partial_{int}C \cap I = \emptyset$, giving (\ref{contour.prop.3}).

It is clear that for all $x \in
\partial_{int}R$ there is $y \in R \cap I$ with $x$ adjacent to $y$; so
(\ref{contour.prop.4}) follows from $\partial_{int} W \subseteq
\partial_{int} R$.

Since $\partial_{int} W \subseteq \cO$, we have $W^\cO \supseteq
\partial_{ext} W^\cE$. If there is $y \in W^\cO$ with $\partial  y \cap W^\cE = \emptyset$,
then the connectivity of $W$ implies that $W=W^\cO$ (and that
$W^\cO$ consists of a single vertex). But $W^\cE$ is non-empty; so
we get the reverse containment $W^\cO \subseteq
\partial_{ext} W^\cE$.

The containment $W^\cE \subseteq \{y \in \cE:\partial  y \subseteq
W^\cO\}$ follows immediately from $W^\cO \supseteq
\partial_{ext} W^\cE$. For the reverse containment,
consider (for a contradiction) $y \in \cE$ with $\partial  y
\subseteq W^\cO$ but $y \not \in W^\cE$. We must have $y \in C$; but
$y$ is not adjacent to anything else in $C$, and $|C| > 1$ (indeed,
$|C|\geq |V|/2 > 1$ since $\grg$ is enveloping), a contradiction
since $C$ is connected. So we have $W^\cE \supseteq \{y \in
\cE:\partial y \subseteq W^\cO\}$. \qed

\medskip

We now return to $T_{L,d}$. Set $\cI_{even}=\{I \in
\cI(T_{L,d}):I~\mbox{even}\}$ and define $\cI_{odd}$ analogously.
The next lemma establishes some of the geometric properties of
$T_{L,d}$ that we will need. Before stating it we need some more
notation.

For $k\geq 1$, we say that $S \subseteq V(T_{L,d})$ is {\em
$k$-clustered} if for every $x, y \in S$ there is a sequence
$x=x_0,\ldots,x_m=y$ of vertices of $S$ such that $d(x_{i-1},x_i)
\leq k$ for all $i=1,\ldots m$, where $d(\cdot,\cdot)$ is the usual
graph distance. Note that $S$ can be partitioned uniquely into
maximal $k$-clustered subsets; we refer to these as the {\em
$k$-components} of $S$.

For a cutset $\gamma$, we define a graph $G_\grg$ as follows. The
vertex set of $G_\grg$ is the set of edges of $T_{L,d}$ that
comprise $\grg$. Declare $e, f \in \gamma$ to be adjacent in
$G_\grg$ if either $e$ and $f$ share exactly one endpoint and if the
coordinate on which the endpoints of $e$ differ is different from
the coordinate on which the endpoints of $f$ differ ({\em i.e.}, $e$
and $f$ are not parallel) or if the endpoints of $e$ and $f$
determine a cycle of length four (a square) in $T_{L,d}$. (This is
equivalent to the following construction, well known in the
statistical physics literature: for $e \in \gamma$, let $e^\star$ be
the dual ($d-1$)-dimensional cube which is orthogonal to $e$ and
bisects it when $T_{L,d}$ is considered as immersed in the continuum
torus. Then declare $e, f \in \gamma$ to be adjacent if $e^\star
\cap f^\star$ is a ($d-2$)-dimensional cube.) We say that a cutset
$\gamma$ is {\em trivial} if $G_\grg$ has only one component.

\begin{lemma} \label{lem-properties.of.contours_specific}
For each $I \in \cI_{even}$ and $\grg \in \Gamma(I)$,
\begin{equation} \label{contour.prop.7}
|\grg| \geq |W|^{1-1/d};
\end{equation}
\begin{equation} \label{contour.prop.8}
\mbox{for large enough $d$, } |\grg| \geq d^{1.9};
\end{equation}
\begin{equation} \label{lem-nontrivial}
\mbox{if $\gamma$ is not trivial then each component of $G_\grg$ has
at least $L^{d-1}$ edges}
\end{equation}
and
\begin{equation} \label{lem-nontrivial2}
\mbox{either $\partial_{int}W$ is $2$-clustered or each of its
$2$-components has size at least $L^{d-1}/2d$}.
\end{equation}
\end{lemma}

\noindent {\em Proof: }For (\ref{contour.prop.7}) and
(\ref{contour.prop.8}) we appeal to an isoperimetric inequality of
Bollob\'as and Leader \cite{BollobasLeader} which states that if $A
\subseteq V(T_{L,d})$ with $|A| \leq L^d/2$, then
$$
|\partial_{ext}A| \geq \min
\left\{2|A|^{1-1/r}rL^{(d/r)-1}:r=1,\ldots,d\right\}.
$$
From this (\ref{contour.prop.7}) follows easily, as does
(\ref{contour.prop.8}) once we observe that $|W| \geq 2d+1$ (since
$W^\cE \neq \emptyset$) and that $|\grg|\geq |\partial_{ext}W|$.

From \cite[Lemma 3]{BorgsChayesFriezeKimTetaliVigodaVu} we have
(\ref{lem-nontrivial}). Finally we turn to (\ref{lem-nontrivial2}).
Let $C_1, \ldots, C_\ell$ be the components of $G_\grg$, and for
each $i$ let $C^\prime_i$ be the vertices of $\partial_{int}W$ which
are endpoints of edges of $C_i$. It is readily checked that each
$C^\prime_i$ is $2$-clustered and that $\partial_{int}W = \cup_i
C^\prime_i$. If $\ell=1$ we therefore have that $\partial_{int}W$ is
$2$-clustered. If $\ell > 1$, we have (by (\ref{lem-nontrivial}))
that each $C_i$ has at least $L^{d-1}$ edges. Since each vertex in
$T_{L,d}$ has degree $2d$, it follows that each $C^\prime_i$ has
size at least $L^{d-1}/2d$. Since the $C^\prime_i$'s are
$2$-clustered, each $2$-component of $\partial_{int} W$ has size at
least $L^{d-1}/2d$, establishing (\ref{lem-nontrivial2}).\qed

\subsection{The main lemma} \label{subsec-results.from.GalvinKahn}

For $c \in {\mathbb N}$ and $v \in V(T_{L,d})$ set
$$
\cW(c,v) = \left\{\grg~:~\mbox{$\grg \in \Gamma(I)$ for some $I \in
\cI_{even}$, $|\grg|=c$, $v \in W^\cE$}\right\}
$$
and set $\cW = \cup_{c,v} \cW(c,v)$. A {\em profile} of a collection
$\{\grg_1, \ldots, \grg_\ell\} \subseteq \cW$ is a vector
$\ul{p}=(c_1 ,v_1, \ldots, c_\ell, v_\ell)$ with $\grg_i \in
\cW(c_i,v_i)$ for all $i$.
%Say that the collection is {\em properly
%nested} if the $\grg_i$'s (viewed as sets of edges) are pairwise
%disjoint and satisfy that for $i < j$ one of $\ant \grg_i \subseteq
%\ant \grg_j$, $\ant \grg_i \cap \ant \grg_j = \emptyset$ holds.
Given a profile vector $\ul{p}$ set
$$
\cI(\ul{p}) = \{I \in \cI_{even}:~\mbox{$\Gamma(I)$ contains a
subset with profile $\ul{p}$}\}.
$$
Our main lemma is the following.
\begin{lemma} \label{lem-volume.bounds.from.gk}
There are constants $c, c', d_0>0$ such that the following holds.
For all even $L\geq 4$, $d \geq d_0$, $\gl$ satisfying
(\ref{inq-bound.on.gl}) and profile vector $\ul{p}$,
\begin{equation} \label{inq-bound.on.inner}
\pi_\gl(\cI(\ul{p})) \leq
\exp\left\{-\frac{c'\gb(\gl)\sum_{i=1}^\ell c_i}{d} \right\},
\end{equation}
where $\gb(\gl)=2\log (1+\gl)-\log(1+2\gl)$.
\end{lemma}

This may be thought of as an extension of the main result of
\cite{GalvinKahn}, which treats only $\ell=1$ and in a slightly less
general setting. We will derive Theorem \ref{thm-main2} from Lemma
\ref{lem-volume.bounds.from.gk} in Section \ref{subsec-proof_of_thm}
before proving the lemma in Section \ref{subsec-proofs}. From here
on we assume that the conditions of Theorem \ref{thm-main2} and
Lemma \ref{lem-volume.bounds.from.gk} are satisfied (with $c$ and
$d_0$ sufficiently large to support our assertions). All constants
implied in $O$ and $\Omega$ statements will be absolute. When it
makes no difference to do otherwise, we assume that all large
numbers are integers. We note for future reference that for $\gl$
satisfying (\ref{inq-bound.on.gl}) we have
\begin{equation} \label{inq-bounds_implied_by_gl}
\frac{\gl}{1+\gl} = \omega\left(\frac{1}{d^{1/4}}\right)
~~~~~\mbox{and}~~~~~\gb(\gl) = \omega\left(\frac{1}{d^{1/2}}\right).
\end{equation}

\subsection{The proof of Theorem \ref{thm-main2}} \label{subsec-proof_of_thm}

We begin with an easy count that dispenses with small balanced
independent sets. Set
$$
\cI_{small} = \left\{I \in \cI_{b,\rho}:\min\{|I^\cE|,|I^\cO|\} \leq
L^d/4d^{1/2}\right\}.
$$
and $\cI_{large} = \cI_{b,\rho} \setminus \cI_{small}$.
\begin{lemma} \label{lem-bounding_small}
$$
\pi_\gl(\cI_{small}) \leq \exp\left\{-\frac{3L^{d-1}}{d^4\log^2
L}\right\}.
$$
\end{lemma}

\medskip

\noindent{\em Proof:} We need a well-known result of Chernoff
\cite{Chernoff} (see also \cite{Bollobas3}, p.11). Let $X_1, \ldots,
X_n$ be i.i.d. Bernoulli random variables with ${\bf P}(X_1=1)=p$.
Then for $k \leq pn$
$$
{\bf P}\left(\sum_{i=1}^n X_i \leq k\right) \leq
2^{nH_p\left(\frac{k}{n}\right)}
$$
where $H_p(x) = x\log_2(p/x)+(1-x)\log_2((1-p)/(1-x))$. Note that
$H_p(x)=H(x)+x\log_2 p +(1-x)\log_2(1-p)$ where $H(x)$ is the usual
binary entropy function. Taking $p=\gl/(1+\gl)$ we see that for a
set $X$ with $|X|=n$ and for $c \leq \gl/(1+\gl)$,
\begin{eqnarray*}
\sum_{A \subseteq X,~|A|\leq cn} \frac{\gl^{|A|}}{(1+\gl)^n} & \leq
& 2^{nH_{\gl/(1+\gl)}(c)} \\
& = & 2^{n\left(H(c)+c\log_2 \frac{\gl}{1+\gl} +(1-c)\log_2
\frac{1}{1+\gl}\right)} \\
& = & 2^{n(H(c)+c\log_2 \gl - \log_2 (1+\gl))}
\end{eqnarray*}
from which it follows that
\begin{equation} \label{inq-weighted.binomial}
\sum_{A \subseteq X,~|A|\leq cn} \gl^{|A|} \leq 2^{n(H(c)+c\log_2
\gl)}.
\end{equation}

Now using $(1+\gl)^{L^d/2}$ as a trivial lower bound on $\sum_{I \in
\cI(T_{L,d})} \gl^{|I|}$ and with the subsequent inequalities
justified below, we have
\begin{eqnarray}
\pi_\gl(\cI_{small}) & \leq &
    2\left(\sum_{A \subseteq \cE,~|A|\leq L^d/4d^{1/2}}
    \gl^{|A|}\right)\left(\sum_{B \subseteq \cO,~|B| \leq
    (1/2d^{1/2}+\rho)L^d/2} \gl^{|B|}\right)(1+\gl)^{-L^d/2} \nonumber \\
& \leq &
    \frac{2\exp_2\left\{\frac{L^d}{2}\left(H\left(\frac{1}{2d^{1/2}}\right) + H\left(\frac{1}{2d^{1/2}}+\rho\right) +
    \left(\frac{1}{d^{1/2}}+\rho\right)\log_2
    \gl\right)\right\}}{(1+\gl)^{L^d/2}}
    \label{using.weighted.binomial} \\
& \leq &
    \exp\left\{-\frac{2L^{d-1}}{d^4\log^2
L}\right\}.
    \label{bounding-small}
\end{eqnarray}
In (\ref{using.weighted.binomial}) we use
(\ref{inq-weighted.binomial}) (legitimate since $1/2d^{1/2} \leq
\gl/(1+\gl)$ and $1/2d^{1/2}+\rho \leq \gl/(1+\gl)$, the former by
(\ref{inq-bounds_implied_by_gl}) and the latter by
(\ref{inq-condition.on.alpha.rho.1})); (\ref{bounding-small})
follows from (\ref{inq-condition.on.alpha.rho}). \qed

\medskip

Set $\cI_{large,~even}=\cI_{large} \cap \cI_{even}$ and define
$\cI_{large,~odd}$ analogously. By Lemma \ref{lem-even_or_odd}
$\cI_{large}=\cI_{large,~even} \cup \cI_{large,~odd}$ and by
symmetry $\pi_\gl(\cI_{large,~even})=\pi_\gl(\cI_{large,~odd})$. In
the presence of Lemma \ref{lem-bounding_small}, Theorem
\ref{thm-main2} reduces to bounding (say)
\begin{equation} \label{inq-remaining}
\pi_\gl(\cI_{large,~even}) \leq
\exp\left\{-\frac{3L^{d-1}}{d^4\log^2 L}\right\}.
\end{equation}
Set $\cI_{large,~even}^{non-trivial} =\{I \in \cI_{large,~even} :
~\mbox{there is $\grg \in \Gamma(I)$ with $|\grg|\geq L^{d-1}$}\}$
and $\cI_{large,~even}^{trivial} = \cI_{large,~even} \setminus
\cI_{large,~even}^{non-trivial}$. With the sum below running over
all vectors $\ul{p}$ of the form $(c,v)$ with $v \in V(T_{L,d})$ and
$c \geq L^{d-1}$, and with the inequalities justified below, we have
\begin{eqnarray}
\pi_\gl(\cI_{large,~even}^{non-trivial}) & \leq & \sum_{\ul{p}}
\pi_\gl(\cI(\ul{p})) \nonumber \\
& \leq &
L^{2d}\exp\left\{-\Omega\left(\frac{L^{d-1}\gb(\gl)}{d}\right)\right\}
\label{large1} \\
& \leq &
\exp\left\{-\Omega\left(\frac{L^{d-1}}{d^{3/2}}\right)\right\}
\label{large2}
\end{eqnarray}
We have used Lemma \ref{lem-volume.bounds.from.gk} in (\ref{large1})
and the factor of $L^{2d}$ is for the choices of $c$ and $v$. In
(\ref{large2}) we have used (\ref{inq-bounds_implied_by_gl}).

For $I \in \cI_{large,~even}^{trivial}$ and $\grg \in \Gamma(I)$ we
have $|\gamma| \geq |\ant \grg|^{1-1/d}$ (by (\ref{contour.prop.7}))
and so
$$
\sum_{\grg \in \Gamma(I)} |\gamma|^{d/(d-1)} \geq \sum_{\grg \in
\Gamma(I)} |\ant \gamma| \geq |I^\cE| \geq L^d/4d^{1/2}.
$$
The second inequality is from Lemma \ref{lem-substantial} and the
third follows since $I \not \in \cI_{small}$.

Set $\Gamma_i(I) = \{\grg \in \Gamma(I):2^{i-1} \leq |\grg| <
2^i\}$. Note that $\Gamma_i(I)$ is empty for $2^i < d^{1.9}$ (recall
(\ref{contour.prop.8})) and for $2^{i-1} > L^{d-1}$ so we may assume
that
\begin{equation} \label{inner.property.1}
1.9 \log d \leq i \leq (d-1)\log L+1.
\end{equation}
Since $\sum_{m=1}^\infty 1/m^2 = \pi^2/6$, there is an $i$ such that
\begin{equation} \label{quad.scale}
\sum_{\grg \in \Gamma_i(I)} |\gamma|^{\frac{d}{d-1}} \geq
\Omega\left(\frac{L^d}{d^{1/2}i^2}\right).
\end{equation}
Choose the smallest such $i$ set $\ell=|\Gamma_i(I)|$. We have
$\sum_{\grg \in \Gamma_i(I)} |\grg| \geq \Omega(\ell 2^i)$ (this
follows from the fact that each $\grg \in \Gamma_i(I)$ satisfies
$|\grg| \geq 2^{i-1}$) and
\begin{equation} \label{inner.property.2}
O\left(\frac{dL^d}{2^i}\right) \geq \ell \geq
\Omega\left(\frac{L^d}{2^{\frac{id}{d-1}}i^2d^{1/2}}\right).
\end{equation}
The first inequality follows from that fact that $\sum_\grg |\grg|
\leq dL^d=|E(T_{L,d})|$; the second follows from (\ref{quad.scale})
and the fact that each $\grg$ has $|\grg|^{d/(d-1)} \leq
2^{di/(d-1)}$. We therefore have $I \in \cI(\ul{p})$ for some
$\ul{p}=(c_1, v_1, \ldots, c_\ell, v_\ell)$ with $\ell$ satisfying
(\ref{inner.property.2}), with
\begin{equation} \label{inner.property.3}
\sum_{j=1}^\ell c_j \geq O(\ell 2^i),
\end{equation}
with
\begin{equation} \label{inner.property.4}
c_j \leq 2^i
\end{equation}
for each $j$ and with $i$ satisfying (\ref{inner.property.1}). With
the sum below running over all profile vectors $\ul{p}$ satisfying
(\ref{inner.property.1}), (\ref{inner.property.2}),
(\ref{inner.property.3}) and (\ref{inner.property.4}) we have
\begin{eqnarray}
\pi_\gl(\cI_{large,~even}^{trivial}) & \leq & \sum_{\ul{p}}
\pi_\gl(\cI(\ul{p})) \nonumber \\
& \leq & d\log L~\max_{i~\mbox{satisfying
(\ref{inner.property.1})}}2^{\ell i}{L^d \choose
\ell}\exp\left\{-\Omega\left(\frac{\ell
2^i\gb(\gl)}{d}\right)\right\}. \label{triv1}
\end{eqnarray}
In (\ref{triv1}) we have used Lemma \ref{lem-volume.bounds.from.gk}.
The factor of $d\log L$ is an upper bound on the number of choices
for $i$; the factor of $2^{\ell i}$ is for the choice of the
$c_j$'s; and the factor ${L^d \choose \ell}$ is for the choice of
the $\ell$ (distinct) $v_j$'s. By (\ref{inner.property.1}), the
second inequality in (\ref{inner.property.2}) and the second
inequality in (\ref{inq-bounds_implied_by_gl}) we have (for $d$
sufficiently large)
\begin{eqnarray*}
2^{\ell i}{L^d \choose \ell} & \leq & 2^{\ell
i}\left(\frac{L^d}{\ell}\right)^\ell \\
& \leq & 2^{\ell i} \left(O\left(2^\frac{id}{d-1}i^2
d^{1/2}\right)\right)^\ell \\
& \leq & 2^{4\ell i} \\
& = & \exp\left\{o\left(\frac{2^i \gb(\gl)}{d}\right)\right\}.
\end{eqnarray*}
Inserting into (\ref{triv1}) we finally get
\begin{eqnarray}
\pi_\gl(\cI_{large,~even}^{trivial}) & \leq & d \log L ~\max_i
\exp\left\{-\Omega\left(\frac{2^i\gb(\gl) \ell}{d}\right)\right\}
\nonumber \\
& \leq & d \log L ~\max_i
\exp\left\{-\Omega\left(\frac{2^i\gb(\gl)L^d}{d
2^{\frac{id}{d-1}}i^2d^{1/2}}\right)\right\} \label{local_2} \\
& \leq & \exp\left\{-\frac{4L^{d-1}}{d^4 \log^2 L}\right\}
\label{local_3}.
\end{eqnarray}
In (\ref{local_2}) we have taken $\ell$ as small as possible, and in
(\ref{local_3}) we have taken $i$ as large as possible and used
(\ref{inq-bounds_implied_by_gl}).

Combining (\ref{local_3}) and (\ref{large2}) we obtain
(\ref{inq-remaining}) and so Theorem \ref{thm-main2}.

\section{Proof of Lemma \ref{lem-volume.bounds.from.gk}}
\label{subsec-proofs}

Our strategy is the following. Let a profile vector
$\ul{p}=(c_1,v_1, \ldots, c_\ell, v_\ell)$ be given. Set
$\ul{p'}=(c_2 ,v_2, \ldots, c_\ell, v_\ell)$. We will show
\begin{equation} \label{induction}
\frac{\pi_\gl(\cI(\ul{p}))}{\pi_\gl(\cI(\ul{p'}))} \leq
\exp\left\{-\Omega\left(\frac{c_1\gb(\gl)}{d}\right)\right\}.
\end{equation}
Then by a telescoping product
$$
\pi_\gl(\cI(\ul{p})) \leq
\frac{\pi_\gl(\cI(\ul{p}))}{\pi_\gl(\cI_{even})} \leq
\exp\left\{-\Omega\left(\frac{\gb(\gl)\sum_{i=1}^\ell c_i}{d}\right)
\right\}
$$
as claimed. To obtain (\ref{induction}) we employ a general strategy
to bound $\pi_\gl(\cS)/\pi_\gl(\cT)$ for $\cS \subseteq \cT
\subseteq \cI(T_{L,d})$ (note that $\cI(\ul{p}) \subseteq
\cI(\ul{p'})$). We define a one-to-many map $\varphi$ from $\cS$ to
$\cT$. We then define a flow $\nu:\cS \times \cT \rightarrow
[0,\infty)$ supported on pairs $(I,J)$ with $J \in \varphi(I)$
satisfying
\begin{equation} \label{eq-flow.out}
\forall I \in \cS, ~~~\sum_{J \in \varphi(I)} \nu(I,J) =1
\end{equation}
and
\begin{equation} \label{eq-flow.in}
\forall J \in \cT, ~~~\sum_{I \in \varphi^{-1}(J)}
\gl^{|I|-|J|}\nu(I,J) \leq M.
\end{equation}
This gives
\begin{eqnarray*}
\sum_{I \in \cS} \gl^{|I|} & = & \sum_{I \in \cS}
\gl^{|I|} \sum_{J \in \varphi(I)} \nu(I,J) \\
& = & \sum_{J \in \cT} \gl^{|J|} \sum_{I \in
\varphi^{-1}(J)} \gl^{|I|-|J|}\nu(I,J) \\
& \leq & M\sum_{J \in \cT} \gl^{|J|}
\end{eqnarray*}
and so $\pi_\gl(\cS)/\pi_\gl(\cT) \leq M$. So our task is to define
$\varphi$ and $\nu$ for $\cS=\cI(\ul{p})$ and $\cT=\cI(\ul{p'})$ for
which (\ref{eq-flow.in}) holds with $M$ given by the right-hand side
of (\ref{induction}).

\medskip

Much of what follows is modified from \cite{GalvinKahn}.  The main
result of \cite{GalvinKahn} has already been described in Section
\ref{sec-intro}. It will be helpful here to describe the main
technical work of that paper. Let $\Lambda_L$ be the box $[-L,L]^d$
in ${\mathbb Z}^d$ with boundary $\partial^\star \Lambda_L =
[-L,L]^d\setminus [-(L-1),L-1]^d$. Write $\cJ$ for the set of
independent sets in $\Lambda_L$ which extend $\partial^\star
\Lambda_L \cap \cO$ and, for a fixed vertex $v_0 \in \Lambda_L \cap
\cE$, write $\cI$ for those $I \in \cJ$ with $v_0 \in I$. The stated
aim of \cite{GalvinKahn} is to show, using a similar strategy to
that described above, that $\pi_\gl(\cI)/\pi_\gl(\cJ) \leq
(1+\gl)^{-2(d-o(1))}$. More specifically, for each $I \in \cI$ let
$\grg'(I)$ be the cutset associated with that component of
$(I^\cE)^+$ that includes $v_0$. For each $w_o, w_e$ write $\cI(w_o,
w_e)$ for those $I \in \cI$ with $|W^\cE|=w_e$ and $|W^\cO|=w_o$,
where $W$ is the subset of $\Lambda_L$ associated with $\grg'(I)$ as
described in Section \ref{subsec-cutsets}. It is shown in
\cite{GalvinKahn} (inequalities (62) and (63) of that paper) that
for $\gl$ satisfying \eqref{inq-bound.on.gl} we have
\begin{equation} \label{from.gk}
\frac{\pi_\gl(\cI(w_o,w_e))}{\pi_\gl(\cJ)} \leq \left\{
    \begin{array}{ll}
          \exp\{-\Omega(\gl^2(w_o-w_e))\} & \mbox{for $\gl < 2$ and} \\
          \gl^{-\Omega(w_o-w_e)} & \mbox{for larger $\gl$}
    \end{array}
    \right.
\end{equation}
from which the stated bound on $\pi_\gl(\cI)/\pi_\gl(\cJ)$ is easily
obtained by a summation. The remainder of this paper is devoted to
an explanation of how the proof of (\ref{from.gk}) needs to be
augmented and modified to obtain our main lemma, and we do not state
the proofs of many of our intermediate lemmas, since they can be
found in the generality we need in \cite{GalvinKahn}. The main
technical issue we have to deal with in moving from (\ref{from.gk})
to Lemma \ref{lem-volume.bounds.from.gk} relates to dealing with
$\grg$ that are non-trivial (in the sense defined before the proof
of Lemma \ref{lem-properties.of.contours_specific}); this is not an
issue in \cite{GalvinKahn} because it is shown there that the
cutsets $\grg'(I)$ described above are always trivial.

\medskip

One technical issue aside, the specification of $\varphi$ is
relatively straightforward. For each $s \in \{\pm 1, \ldots, \pm
d\}$, define $\sigma_s$, the {\em shift in direction $s$}, by
$\sigma_s(x)=x+e_s$, where $e_s$ is the $s$th standard basis vector
if $s>0$ and $e_s=-e_{-s}$ if $s<0$. For $X \subseteq V(T_{L,d})$,
write $\sigma_s(X)$ for $\{\sigma_s(x):x \in X\}$. For a cutset
$\grg \in \cW$ set $W^s = \{x \in
\partial_{int}W:\sigma_s^{-1}(x) \not \in W\}$. We will obtain
$\varphi(I)$ by shifting $I$ inside $W$ in a certain direction $s$
and adding arbitrary subsets of $W^s$ to the result, where $W$ is
associated with a cutset $\grg \in \Gamma(I) \cap \cW(c_1,v_1)$. The
success of this process depends on the fact that $I$ is disjoint
from the vertex set of $\grg$. We now formalize this.

\begin{lemma} \label{lem-shift_works}
Let $I \in \cI(\ul{p})$ be given. Let $\grg \in \Gamma(I)$ be such
that $|\grg|=c_1$ and $v_1 \in W^\cE$ where $W=\ant \grg$. For any
choice of $s$, it holds that
$$
I_0 := (I\setminus W) \cup \sigma_s(I \cap W) ~~\mbox{is in}~~
\cI(\ul{p'})
$$
and has the same size as $I$. Moreover, the sets $I_0$ and $W^s$ are
mutually disjoint and
$$
I_0 \cup W^s \in \cI(\ul{p'}).
$$
\end{lemma}

\noindent {\em Proof: }That $I_0 \cup W^s$ is an independent set and
that $I_0$ is the same size as $I$ is the content of
\cite[Proposition 2.12]{GalvinKahn}. Because $\ant \grg$ is disjoint
from the interiors of the remaining cutsets and the shift operation
that creates $I_0 \cup W^s$ only modifies $I$ inside $W$ it follows
that $I_0, I_0 \cup W^s \in \cI(\ul{p'})$. \qed

\medskip

For $I \in \cI(\ul{p})$ we define
$$
\varphi(I)=\{I_0 \cup S:S \subseteq W^s\}
$$
for a certain $s$ to be chosen presently. In light of Lemma
\ref{lem-shift_works}, $\varphi(I) \subseteq \cI(\ul{p'})$
regardless of this choice.

\medskip

To define $\nu$ and $s$ we employ the notion of approximation also
used in \cite{GalvinKahn} and introduced by Sapozhenko in
\cite{Sap}. For $\grg \in \cW$ we say that $A \subseteq V(T_{L,d})$
is an {\em approximation} of $\grg$ if
\begin{equation} \label{approxpair0}
A^\cE \supseteq W^\cE ~~~\mbox{and}~~~A^\cO \subseteq W^\cO,
\end{equation}
\begin{equation} \label{approxpair1}
d_{A^\cO}(x) \geq 2d-\sqrt{d} ~~\mbox{for all $x \in A^\cE$}
\end{equation}
and
\begin{equation} \label{approxpair2}
d_{\cE \setminus A^\cE}(x) \geq 2d-\sqrt{d} ~~\mbox{for all $y \in
\cO \setminus A^\cO$},
\end{equation}
where $d_X(x)=|\partial x \cap X|$. Note that since $W_\cO =
\partial W_\cE$, $W$ is an approximation of
$\grg$.

To motivate the definition of approximation, note that by
(\ref{contour.prop.5}) if $u$ is in $W^\cE$ then all of its
neighbors are in $W^\cO$, and if $u'$ is in $\cO \setminus W^\cO$
then all of its neighbors are in $\cE \setminus W^\cE$. If we think
of $A^\cE$ as approximate-$W^\cE$ and $A^\cO$ as
approximate-$W^\cO$, (\ref{approxpair1}) says that if $u \in \cE$ is
in approximate-$W^\cE$ then almost all of its neighbors are in
approximate-$W^\cO$ while (\ref{approxpair2}) says that if $u'\in
\cO$ is not in approximate-$W^\cO$ then almost all of its neighbors
are not in approximate-$W^\cE$.

\medskip

Before stating our main approximation lemma, which is a slight
modification of \cite[Lemma 2.18]{GalvinKahn}, it will be convenient
to further refine our partition of cutsets. To this end set
$$
\cW(w_e,w_o,v) = \{\grg: \mbox{$\grg \in \Gamma(I)$ for some $I \in
\cI_{even}$, $|W^\cO|=w_o$, $|W^\cE|=w_e$, $v \in W^\cE$}\}.
$$
Note that (by (\ref{contour.prop.5}))
$$
|\grg| = |\nabla(W)|=2d(|W^\cO|-|W^\cE|)
$$
so $\cW(w_e,w_o,v) \subseteq \cW((w_o-w_e)/2d, v)$.

\begin{lemma} \label{lem-k.comp.approx}
For each $w_e$, $w_o$ and $v$ there is a family $\cA(w_e,w_o,v)$
satisfying
$$
|\cA(w_e,w_o,v)| \leq
\exp\left\{O\left((w_o-w_e)d^{-\frac{1}{2}}\log^{\frac{3}{2}}d\right)\right\}
$$
and a map $\Pi:\cW(w_e,w_o,v) \rightarrow \cA(w_e,w_o,v)$ such that
for each $\grg \in \cW(w_e,w_o,v)$, $\Pi(\grg)$ is an approximation
for $\grg$.
\end{lemma}

The proof of this lemma is deferred to Section
\ref{subsec-proof.of.approx.lemma}. Our bound on the number of
approximate cutsets with parameters $w_e$, $w_o$ and $v$ is much
smaller than any bound we are able to obtain on the number of
cutsets with the same set of parameters. This is where we make the
entropy gain discussed in Section \ref{sec-overview}.

\medskip

We are now in a position to define $\nu$ and $s$. Our plan for each
fixed $J \in \cI(\ul{p'})$ is to fix $w_e, w_o$ and $A \in \cW(w_e,
w_o)$ and to consider the contribution to the sum in
(\ref{eq-flow.in}) from those $I \in \varphi^{-1}(J)$ with
$\Pi(\grg(I))=A$. We will try to define $\nu$ in such a way that
each of these individual contributions to (\ref{eq-flow.in}) is
small; to succeed in this endeavour we must first choose $s$ with
care. To this end, given $\grg \in \cW(w_e,w_o,v)$, set
$$
Q^\cE = A^\cE \cap \partial_{ext}(\cO \setminus A^\cO)
~~~~~\mbox{and}~~~~~Q^\cO = (\cO \setminus A^\cO) \cap
\partial_{ext} A^\cE,
$$
where $A = \Pi(\grg)$ in the map guaranteed by Lemma
\ref{lem-k.comp.approx}. To motivate the introduction of $Q^\cE$ and
$Q^\cO$, note that for $\grg \in \Pi^{-1}(A)$ we have
\begin{eqnarray*}
A^\cE\setminus Q^\cE & \subseteq & W^\cE \\
\cE \setminus A^\cE & \subseteq & \cE \setminus W^\cE \\
A^\cO & \subseteq & W^\cO
\end{eqnarray*}
and
\begin{eqnarray*}
\cO \setminus (A^\cO \cup Q^\cO) & \subseteq & \cO \setminus W^ \cO
\end{eqnarray*}
(all using \eqref{contour.prop.5} and \eqref{approxpair0}). It
follows that for each $\grg \in \Pi^{-1}(A)$, $Q^\cE \cup Q^\cO$
contains all of the vertices whose location in the partition
$T_{L,d} = W \cup \overline{W}$ is as yet unknown.

%\begin{lemma} \label{lem-ambiguous}
%Fix $A \in \cA(w_e,w_o,v)$. If $u \in V(T_{L,d})$ has the property
%that there exists $\grg_1,\grg_2 \in \Pi^{-1}(A)$ with $u \in W_1$
%and $u \not \in W_2$ then $x \in Q^\cE \cup Q^\cO$.
%\end{lemma}
%
%\noindent {\em Proof: }\cite[discussion after statement of Lemma
%2.18]{GalvinKahn}. \qed
%
%\medskip
%
%Motivated by Lemma \ref{lem-ambiguous} we refer to $Q^\cE \cup
%Q^\cO$ as the {\em ambiguous} vertices with respect to $A$.

\begin{lemma} \label{lem-lattice_direction}
For $\grg \in \cW(w_e,w_o,v)$, there is an $s \in \{\pm 1, \ldots,
\pm d\}$ such that both of
$$
|W^s| \geq .8(w_o-w_e)~~~~~\mbox{and}~~~~~ |\sigma_s(Q^\cE) \cap
Q^\cO| \leq \frac{5|W^s|}{\sqrt{d}}
$$
hold.
\end{lemma}

\noindent {\em Proof: }\cite[(49) and (50)]{GalvinKahn}. \qed

\medskip

We choose the smallest such $s$ to be the lattice direction
associated with $\grg$. Note that $s$ depends on $\grg$ but not on
$I$.

Now for each $I \in \cI(\ul{p})$ let $\grg \in \Gamma(I)$ be a
particular cutset with $\grg \in \cW(c_1, v_1)$. Let $\varphi(I)$ be
as defined before, with $s$ as specified by Lemma
\ref{lem-lattice_direction}. Define
$$
C = W^s \cap A^\cO \cap \gs_s(Q^\cE)
$$
and
$$
D = W^s \setminus C,
$$
and for each $J \in \varphi(I)$ set
$$
\nu(I,J) = \gl^{|J \cap W^s|} \left(\frac{\gl}{(1+\gl)^2}\right)^{|C
\cap J|} \left(\frac{1+2\gl}{(1+\gl)^2}\right)^{|C \setminus
J|}\left(\frac{1}{1+\gl}\right)^{|D|}.
$$
Note that for $I \in \varphi^{-1}(J)$, $\nu(I,J)$ depends on $W$ but
not on $I$ itself.

Noting that $C \cup D$ partitions $W$ we have
\begin{eqnarray}
\sum_{J \in \varphi(I)} \nu(I,J) & = & \sum_{A\subseteq C,~B
\subseteq D} \gl^{|A|+|B|}\left(\frac{\gl}{(1+\gl)^2}\right)^{|A|}
\left(\frac{1+2\gl}{(1+\gl)^2}\right)^{|C|-|A|}\left(\frac{1}{1+\gl}\right)^{|D|}
\nonumber \\
& = & \sum_{B\subseteq D}\frac{\gl^{|B|}}{(1+\gl)^{|D|}}\sum_{A
\subseteq
C}\left(\frac{\gl^2}{1+2\gl}\right)^{|A|}\left(\frac{1+2\gl}{(1+\gl)^2}\right)^{|C|}
\nonumber \\
& = & \frac{(1+\gl)^{|D|}}{(1+\gl)^{|D|}}
\left(\frac{1+2\gl+\gl^2}{1+2\gl}\right)^{|C|}\left(\frac{1+2\gl}{(1+\gl)^2}\right)^{|C|}
 \nonumber
\\
& = & 1, \nonumber
\end{eqnarray}
so $\nu$ satisfies (\ref{eq-flow.out}). To obtain
(\ref{inq-bound.on.inner}) we must establish (\ref{eq-flow.in}) with
$M$ given by the right-hand side of (\ref{induction}).

Fix $w_e$, $w_o$ such that $2d(w_o-w_e)=c_1$. Fix $A \in
\cA(w_e,w_o,v_1)$ and $s \in \{\pm 1, \ldots, \pm d\}$. For $I$ with
$\grg(I) \in \cW(w_e,w_o, v_1)$ write $I \sim_s A$ if it holds that
$\Pi(\grg)=A$ and $s(I)=s$. The next lemma, which bounds the
contribution to the sum in (\ref{eq-flow.in}) from those $I \in
\varphi^{-1}(J)$ with $I \sim_s A$, is the heart of the whole proof,
and perhaps the principal inequality of \cite{GalvinKahn}. We
extract it directly from \cite{GalvinKahn}; although the setting
here is slightly different, the proof is identical to the equivalent
statement in \cite{GalvinKahn}.
\begin{lemma} \label{lem-gk_black_box}
For $J \in \cI(\ul{p'})$,
$$
\sum \left\{\gl^{|I|-|J|}\nu(I,J):I \sim_s A,~I \in
\varphi^{-1}(J)\right\} \leq
\left(\frac{\sqrt{1+2\gl}}{1+\gl}\right)^{w_o-w_e}.
$$
\end{lemma}

\noindent {\em Proof: }\cite[Section 2.12]{GalvinKahn}.
\qed

\medskip

We are now only a short step away from (\ref{inq-bound.on.inner}).
%
%In Section \ref{subsec-outline_of_(53)} we will give an outline of
%the key steps of the proof of Lemma \ref{lem-gk_black_box} (these
%are not spelled out explicitly in \cite{GalvinKahn}), but first we
%derive (\ref{inq-bound.on.inner}).
With the steps justified below we have that for each $J \in
\cI(\ul{p'})$
\begin{eqnarray}
\sum_{I \in \varphi^{-1}(J)} \gl^{|I|-|J|}\nu(I,J) & \leq &
\sum_{w_e,w_o}\sum_{s, A} \sum \left\{\gl^{|I|-|J|}\nu(I,J):I \sim_s
A,~I \in
\varphi^{-1}(J)\right\} \nonumber \\
& \leq & 2d c_1^{\frac{2d}{d-1}}|\cA(w_e, w_o,
v_1)|\left(\frac{\sqrt{1+2\gl}}{1+\gl}\right)^{\frac{c_1}{2d}}
\label{overcount} \\
& \leq & 2d
c_1^{\frac{2d}{d-1}}\exp\left\{-\Omega\left(\frac{c_1\gb(\gl)}{d}\right)\right\}
\label{overcount2} \\
& \leq &
\exp\left\{-\Omega\left(\frac{c_1\gb(\gl)}{d}\right)\right\}
\label{overcount3}
\end{eqnarray}
completing the proof of (\ref{eq-flow.in}). In (\ref{overcount}), we
note that there are $|\cA(w_e, w_o, v_1)|$ choices for the
approximation $A$, $2d$ choices for $s$ and $c_1^{d/(d-1)}$ choices
for each of $w_e$, $w_o$ (this is because $c_1\geq
(w_e+w_o)^{1-1/d}$ by (\ref{contour.prop.7})), and we apply Lemma
\ref{lem-gk_black_box} to bound the summand. In (\ref{overcount2})
use Lemma \ref{lem-k.comp.approx} and the fact that for any $c>0$ we
can choose $c^\prime>0$ such that whenever $\gl
> c^\prime d^{-1/4}\log^{3/4}d$ and $d=d(c)$ is sufficiently large we have
$$
\exp\left\{cd^{-\frac{1}{2}}\log^{\frac{3}{2}}d\right\}\frac{\sqrt{1+2\gl}}{1+\gl}
\leq \exp\left\{-\frac{\gb(\gl)}{4}\right\}.
$$
Finally in (\ref{overcount3}) we use $c_1 \geq d^{1.9}$ (by
(\ref{contour.prop.8})) and the second inequality in
(\ref{inq-bounds_implied_by_gl}) to bound
$2dc_1^{2d/(d-1)}=\exp\{o(c_1 \gb(\gl)/d)\}$.

\subsection{Proof of Lemma \ref{lem-k.comp.approx}}
\label{subsec-proof.of.approx.lemma}

We obtain Lemma \ref{lem-k.comp.approx} by combining a sequence of
lemmas. Lemma \ref{lem-approx_step.1}, which we extract directly
from \cite{GalvinKahn}, establishes the existence for each $\grg$ of
a very small set of vertices nearby to $\grg$ whose neighbourhood
can be thought of as a coarse approximation to $\grg$. (We will
elaborate on this after the statement of the lemma.) Lemma
\ref{lem-approx_step.2} shows that there is a small collection of
these coarse approximations such that every $\grg \in
\cW(w_e,w_o,v)$ is approximated by one of the collection. Our proof
of this lemma for $\grg$ trivial is from \cite{GalvinKahn}, but we
need to add a new ingredient to deal with non-trivial $\grg$.
Finally Lemma \ref{lem-approx_step.3}, which we extract directly
from \cite{GalvinKahn}, turns the coarse approximations of Lemma
\ref{lem-approx_step.2} into the more refined approximations of
Lemma \ref{lem-k.comp.approx} without increasing the number of
approximations too much.

Given $\grg  \in \cW(w_e,w_o,v)$ set
$$
\partial_{int}'W = \{x \in \partial_{int}W : d_{W^\cE}(x) \leq d\}
~~~\mbox{and}~~~
\partial_{int}'C = \{x \in \partial_{int}C : d_{C^\cO}(x) \leq d\}.
$$
(Recall that $d_X(x)=|\partial x \cap X|$.)
\begin{lemma} \label{lem-approx_step.1}
For each $\grg \in \cW(w_e,w_o,v)$ there is a $U$ with the following
properties.
\begin{equation} \label{prop-u.-1}
U \subseteq N(\partial_{int}'W \cup \partial_{int}'C)
\end{equation}
\begin{equation} \label{prop-u.0}
N(U) \supseteq \partial_{int}'W \cup \partial_{int}'C
\end{equation}
and
\begin{equation} \label{prop-u.3}
|U| \leq O\left((w_o-w_e)\sqrt{\frac{\log 2d}{2d}}\right)
\end{equation}
where $N(X)=\cup_{x \in X} \partial x$.
\end{lemma}
To motivate Lemma \ref{lem-approx_step.1}, let us point out that in
\cite[(34)]{GalvinKahn} it is observed that for $U$ satisfying
(\ref{prop-u.-1}) and (\ref{prop-u.0}) the removal of $N(U)$ from
$V(T_{L,d})$ separates $W$ from $C$. $U$ may therefore be thought of
a coarse approximation to $\grg$: removing $U$ and its neighbourhood
achieves the same effect as removing $\grg$. However, $U$ is very
much smaller than $\grg$ ($\grg$ has $2d(w_o-w_e)$ edges). By
focusing on specifying $U$ instead of $\grg$, we lose some
information, but we gain because fewer choices have to be made to
specify $U$. The engine driving the proof of Lemma
\ref{lem-volume.bounds.from.gk} is the fact that the gain far
outweighs the loss. Lemma \ref{lem-approx_step.1} is \cite[Lemma
2.15]{GalvinKahn} and we omit the proof.

\begin{lemma} \label{lem-approx_step.2}
For each $w_e, w_o$ and $v$ there is a family $\cU(w_e,w_o,v)$
satisfying
$$
|\cU(w_e,w_o,v)| \leq
\exp\left\{O\left((w_o-w_e)d^{-\frac{1}{2}}\log^{\frac{3}{2}}d\right)\right\}
$$
and a map $\Pi^\cU:\cW(w_e,w_o,v) \rightarrow \cU(w_e,w_o,v)$ such
that for each $\grg \in \cW(w_e,w_o,v)$, $\Pi^\cU(\grg)$ satisfies
(\ref{prop-u.-1}), (\ref{prop-u.0}) and (\ref{prop-u.3}).
\end{lemma}

\noindent {\em Proof: }It is observed in \cite[paragraph after
(35)]{GalvinKahn} that for $U$ satisfying (\ref{prop-u.-1}) and
(\ref{prop-u.0}) we have
\begin{equation} \label{prop-u.1}
\mbox{for all $x \in \partial_{int}W$, $d(x,U) \leq 2$}
\end{equation}
and
\begin{equation} \label{prop-u.2}
\mbox{for all $y \in U$, $d(y,\partial_{int}W) \leq 2$.}
\end{equation}
Let $U$ satisfy (\ref{prop-u.1}), (\ref{prop-u.2}) and
(\ref{prop-u.3}) for some $\grg \in \cW(w_e,w_o,v)$ and let $W_1,
\ldots, W_k$ be the $2$-components of $\partial_{int}W$. For each
$j=1,\ldots, k$ let
$$
U_j = \{y \in U:d(y,x) \leq 2~~\mbox{for some $x \in W_j$}\}.
$$
We claim that each $U_j$ is $6$-clustered. To see this, fix $u, v
\in W_j$ and take $x_u \in W_j$ at distance at most $2$ from $u$ and
$x_v \in W_j$ at distance at most $2$ from $v$. Let $x_u=x_0,
\ldots, x_\ell=x_v$ be a sequence of vertices in $W_j$ with
$d(x_{i-1},x_i) \leq 2$ for each $i$. For $i=1, \ldots, \ell-1$,
take $u_i \in W_j$ with $d(u_i,x_i) \leq 2$. Then the sequence
$u=u_0,u_1, \ldots, u_{\ell-1},u_\ell=v$ has the property that
$d(u_{i-1},u_i) \leq 6$ for each $i$, establishing the claim.

To bound the number of possibilities for $U$ we first consider the
case $2d(w_o-w_e) \leq L^{d-1}$. In this case, by
(\ref{lem-nontrivial2}), all $\grg$ under consideration are trivial
(in the sense defined before the statement of Lemma
\ref{lem-properties.of.contours_general}) and $k=1$.

We show that there is a small (size $O(w_od^2)$) set of vertices
meeting all possible $U$'s in this case. Fix a linear ordering $\ll$
of $\cO$ satisfying
$$
d(v,y_1) < d(v,y_2)~~~\Longrightarrow~~~y_1 \ll y_2,
$$
and let $T$ be the initial segment of $\ll$ of size $w_o$. We claim
that $T \cap \partial_{int}W \neq \emptyset$. If $T=W^\cO$, this is
clear; if not, consider a shortest $y-v$ path in $T_{L,d}$ for some
$y \in T \setminus W^\cO$. This path intersects $W^\cO$ (since
$\partial v\subseteq W^\cO$). Let $y'$ be the largest (with respect
to $\ll$) vertex of $W^\cO$ on the path; then $y' \in
\partial_{int}W \cap T$, establishing our claim. There are at most
$w_o$ possibilities for $y' \in \partial_{int}W \cap T$, so at most
$O(w_od^2)$ possibilities for a vertex $x'$ with $d(x',y')\leq 2$;
and by (\ref{prop-u.1}) $U$ must contain such an $x'$.

In this case we may take $\cU(w_e,w_o,v)$ to be the collection of
all $6$-connected subsets of $V(T_{L,d})$ of size at most
$O((w_o-w_e)\sqrt{\log 2d / 2d})$ containing one of the $O(w_od^2)$
vertices described in the last paragraph. Using the fact that in any
graph with maximum degree $\Delta$ the number of connected, induced
subgraphs of order $n$ containing a fixed vertex is at most
$(e\Delta)^n$ (see, {\em e.g.}, \cite[Lemma 2.1]{GalvinKahn}) we
infer that
\begin{eqnarray}
|\cU(w_e,w_o,v)| & \leq &
O(w_od^2)(d^7)^{O\left((w_o-w_e)\sqrt{\frac{\log
2d}{2d}}\right)} \label{mod1} \\
& \leq &
\exp\left\{O\left((w_o-w_e)d^{-\frac{1}{2}}\log^{\frac{3}{2}}d\right)\right\},
\label{mod3}
\end{eqnarray}
as required. The factor of $O(w_od^2)$ in (\ref{mod1}) accounts for
the choice of a fixed vertex in $U$; the exponent
$O((w_o-w_e)\sqrt{\log 2d/2d})$ is from (\ref{prop-u.3}); and the
$d^7$ accounts for the fact that $U$ is connected in a graph with
maximum degree at most $65d^6$. In (\ref{mod3}) we use
(\ref{contour.prop.7}) to bound $2d(w_o-w_e) \geq (w_o+w_e)^{1-1/d}
\geq w_o^{3/4}$ and so (since $w_o \geq 2d$) $\log(w_od^2) =
o((w_o-w_e)d^{-1/2}\log^{3/2}d)$.

In the case where $2d(w_o-w_e) > L^{d-1}$, by
(\ref{lem-nontrivial2}) each of the components of $\grg$ has at
least $L^{d-1}$ edges, so $\grg$ has at most $dL^d/L^{d-1}=dL$
components and $U$ at most $dL$ $6$-components. In this case we may
take $\cU(w_e,w_o,v)$ to be the collection of all subsets of
$V(T_{L,d})$ of size at most $O((w_o-w_e)\sqrt{\log 2d / 2d})$
containing at most $dL$ $6$-components. As in the previous case we
have
\begin{eqnarray}
|\cU(w_e,w_o,v)| & \leq &
(L^d)^{dL}(d^7)^{O\left((w_o-w_e)\sqrt{\frac{\log 2d}{2d}}\right)}
\sum_{j=1}^{dL} {O\left(((w_o-w_e)\sqrt{\frac{\log
2d}{2d}}\right)+j-1 \choose j-1}\nonumber \\
& \leq &
\exp\left\{O\left((w_o-w_e)d^{-\frac{1}{2}}\log^{\frac{3}{2}}d\right)\right\},
\label{mod5}
\end{eqnarray}
as required, the extra factors in the first inequality accounting
for the choice of a fixed vertex in each of the at most $dL$
$6$-components and of the sizes of each of the $6$-components. To
obtain (\ref{mod5}) we use $w_o \leq L^d$ to bound
$(L^d)^{dL}\sum_{j=1}^{dL} {O((w_o-w_e)\sqrt{\log 2d/2d})+j-1
\choose j-1} \leq 2^{O(d^2L\log L)}$ and $2d(w_o-w_e)\geq L^{d-1}$
to bound $d^2L\log L = o((w_o-w_e)d^{-1/2}\log^{3/2}d)$. \qed

\medskip

The next lemma turns $\cU(w_e,w_o,v)$ into the collection of
approximations postulated in Lemma \ref{lem-k.comp.approx}. It is a
straightforward combination of \cite[Lemmas 2.16, 2.17,
2.18]{GalvinKahn}, and we omit the proof. Combining Lemmas
\ref{lem-approx_step.2} and \ref{lem-approx_step.3} we obtain Lemma
\ref{lem-k.comp.approx}.

\begin{lemma} \label{lem-approx_step.3}
For each $U \in \cU(w_e,w_o,v)$ there is a family $\cV(w_e,w_o,v)$
satisfying
$$
|\cV(w_e,w_o,v)| \leq
\exp\left\{O\left((w_o-w_e)d^{-\frac{1}{2}}\log^{\frac{3}{2}}d\right)\right\}
$$
and a map $\Pi^\cV:\cU(w_e,w_o,v) \rightarrow \cV(w_e,w_o,v)$ such
that for each $\grg \in \cW(w_e,w_o,v)$ and $U \in \cU(w_e,w_o,v)$
with $\Pi^\cU(\grg)=U$, $\Pi^\cV(U)$ is an approximation of $\grg$.
\end{lemma}

\medskip

\noindent {\bf Acknowledgment: }We thank Dana Randall for numerous
helpful discussions.


\begin{thebibliography}{99}

\bibitem{BergSteif}
J. van den Berg and J.E. Steif, Percolation and the hard-core
lattice gas model, {\em Stoch. Proc. Appl.} {\bf 49} (1994),
179--197.

\bibitem{Bollobas}
B. Bollob\'as, {\em Modern Graph Theory}, Springer, New York,
1998.

\bibitem{Bollobas3} B. Bollob\'as, {\em Random Graphs,} Cambridge University Press,
Cambridge, 2001.

\bibitem{BollobasLeader} B. Bollob\'as and I. Leader, Edge-isoperimetric
inequalities in the grid, {\em Combinatorica} {\bf 11} (1991),
299--314.

\bibitem{BorgsChayesFriezeKimTetaliVigodaVu}
C. Borgs, J. Chayes, A. Frieze, J.H. Kim, P. Tetali, E. Vigoda, V.
Vu, Torpid Mixing of some Monte Carlo Markov Chain algorithms in
Statistical Physics, {\em Proc. of the IEEE FOCS '99}, 218--229.

\bibitem{Chernoff} H. Chernoff, A measure of asymptotic efficiency for tests of
a hypothesis based on the sum of observations, {\em Ann. Math.
Statistics} {\bf 23} (1952), 493--507.

\bibitem{Diestel}
R. Diestel, {\em Graph Theory}, Springer, New York, 2000.

\bibitem{Dobrushin} R.L. Dobrushin,
The problem of uniqueness of a Gibbs random field and the problem of
phase transition, {\em Functional Anal. Appl.} {\bf 2} (1968),
302--312.

\bibitem{Dobrushin2} R.L. Dobrushin,
An investigation of Gibbs states for three-dimensional lattice
systems, {\em Teor. Verojatnost. i Primenen} {\bf 18} (1973),
261--279. (Russian with English summary; translation in Theor. Prob.
Appl. 18 (1974), 253--271.)

\bibitem{DyerFriezeJerrum}
M. Dyer, A. Frieze and M. Jerrum, On counting independent sets in
sparse graphs, {\em SIAM J. Comp.} {\bf 31} (2002), 1527--1541.

\bibitem{DyerGreenhill}
M. Dyer and C. Greenhill, On Markov Chains for independent sets,
{\em Journal of Algorithms} {\bf 35} (2000), 17--49.

\bibitem{GalvinKahn}
D. Galvin and J. Kahn, On phase transition in the hard-core model on
${\mathbb Z}^d$, {\em Comb. Prob. Comp.} {\bf 13} (2004), 137--164.

\bibitem{GalvinTetali}
D. Galvin and P. Tetali, Slow mixing of Glauber dynamics for the
hard-core model on regular bipartite graphs, {\em Random Structures
and Algorithms} {\bf 28} (2006), 427--443.

\bibitem{Grimmett}
G. Grimmett, {\em Percolation}, Springer-Verlag, Berlin, 1999.

\bibitem{JerrumSinclair}
M. Jerrum and A. Sinclair, Conductance and the rapid mixing property
for Markov chains: the approximation of the permanent resolved, {\em
Proc. ACM STOC '88}, 235--243.

\bibitem{Kelly}
F. Kelly, Loss networks, {\em Ann. App. Prob.} {\bf 1} (1991)
319--378.

\bibitem{LebowitzMazel} J. Lebowitz and A. Mazel, Improved Peierl's
argument for higher dimensional Ising models, {\em J. Stat. Physics}
{\bf 90} (1998), 1051--1059.

\bibitem{LubyVigoda}
M. Luby and E. Vigoda, Fast convergence of the Glauber dynamics for
sampling independent sets, {\em Random Structures and Algorithms}
{\bf 15} (1999), 229--241.

\bibitem{MontenegroTetali}
R. Montenegro and P. Tetali, Mathematical aspects of mixing times in
Markov chains, {\em Found. Trends Theor. Comput. Sci.} {\bf 1} no. 3
(2006).

\bibitem{RandallMixing}
D. Randall, Mixing, {\em Proc. of the IEEE FOCS '03}, 4--15.

\bibitem{Sap}
A. A. Sapozhenko, On the number of connected subsets with given
cardinality of the boundary in bipartite graphs, {\em Metody
Diskret. Analiz.} {\bf 45} (1987), 42--70.  (Russian)

\bibitem{Weitz}
D. Weitz, Counting independent sets up to the tree treshold, {\em
Proc. ACM STOC '06}, 140--149.

\end{thebibliography}
\end{document}